\documentclass{article}
\usepackage{amssymb}
\newtheorem{theorem}{Theorem}[section]
\newtheorem{proposition}{Proposition}[section]

\newtheorem{remark}{Remark}[section]

\begin{document}

\title{A sharp estimate and change on the dimension of the attractor for Allen-Cahn equations\footnote{Keywords and Phrases: Allen-Cahn equation, singular potential, Hardy inequality, attractors, Hausdorff dimension, fractal dimension. 
AMS Subject Classification: 35K57, 35B40, 35B41, 37L30.\newline E-mail: karan@aegean.gr (N. I. Karachalios) \& nzogr@science.tuc.gr (N. B. Zographopoulos)}}
\date{}
\author{Nikos. I. Karachalios, \\
{\it Department of Mathematics},\\ 
{\it University of the Aegean},\\
{\it Karlovassi, 83200 Samos, Greece}\\
\\
Nikos B. Zographopoulos,\\
{\it Department of Mathematics, },\\
{\it National Technical University of Crete},\\
{\it Chania, 73100 Crete, Greece}\\
\\
}
\maketitle
\pagestyle{myheadings} \thispagestyle{plain} \markboth{N. B.
Zographopoulos}{}
\maketitle
\begin{abstract}
We consider the semilinear reaction diffusion equation $\partial_t\phi-\nu\Delta\phi-V(x)\phi+f(\phi)=0$, $\nu>0$ in a bounded domain $\Omega\subset\mathbb{R}^N$. We assume the standard Allen-Cahn-type nonlinearity, while the potential $V$ is either the inverse square potential $V(x)=\delta\,|x|^{-2}$ or the  borderline potential $V(x)=\delta\,\mathrm{dist}(x,\partial\Omega)^{-2}$, $\delta\geq 0$ (thus including the classical Allen-Cahn equation as a special case when $\delta=0$). In the subcritical cases $\delta=0$, $N\geq 1$ and $0<\mu:=\frac{\delta}{\nu}<\mu^*$, $N\geq 3$ (where $\mu^*$ is the optimal constant of Hardy and Hardy-type inequalities), we present a new  estimate on the dimension of the global attractor. This estimate comes out by an improved lower bound for sums of eigenvalues of the Laplacian by A. D. Melas (Proc. Amer. Math. Soc. \textbf{131} (2003), 631-636). The estimate is sharp, revealing the existence of (an explicitly given) threshold value for the ratio of the volume to the moment of inertia of $\Omega$ on which the dimension of the attractor may considerably change. Consideration is also given on the finite dimensionality  of the global attractor in the critical case $\mu=\mu^*$.
\end{abstract}
\section{Introduction}
Let $\Omega\subset\mathbb{R}^N$, $N\geq 1$ be a bounded open domain with boundary $\partial\Omega$ and consider the eigenvalues $0<\lambda_1(\Omega)\leq\lambda_2(\Omega)\ldots\leq\lambda_m(\Omega)\leq\ldots$ (counting multiplicity) of the Dirichlet Laplacian
\begin{eqnarray}
\label{D1}
-\Delta u&=&\lambda u,\;\;\mbox{in}\;\;\Omega,\\
u&=&0,\;\;\mbox{on}\;\;\partial\Omega.\nonumber
\end{eqnarray}
In \cite[Theorem 1, pg. 312]{LiYau83} it was proved that for any $m\geq 1$,
\begin{eqnarray}
\label{LiYau1}
\sum_{i=1}^m\lambda_{i}(\Omega)\geq\frac{NC_N}{N+2}\mu_N(\Omega)^{-\frac{2}{N}}m^{\frac{N+2}{N}},\;\;C_N=(2\pi)^2\omega_N^{-2/N}.
\end{eqnarray}
Here $\omega_N$ denotes the volume of the unit ball in $\mathbb{R}^N$ and $\mu_N(\Omega)$ denotes the $N$-dimensional volume of $\Omega$. The lower bound (\ref{LiYau1}) on the sums of the eigenvalues is sharp in view of H. Weyl's asymptotic formula
\begin{eqnarray*}
\lambda_m(\Omega)\sim C_N\left(\frac{m}{\mu_N(\Omega)}\right)^{\frac{2}{N}},\;\;\mbox{as}\;\;m\rightarrow\infty.
\end{eqnarray*}
In \cite[Theorem 1, pg. 632 \& pg. 635]{Melas2002}, the following improvement of (\ref{LiYau1}) is given:  Denote by $I(\Omega)$ the ``moment of inertia'' of $\Omega$, defined as
\begin{eqnarray*}
I(\Omega)=\min_{\alpha\in\mathbb{R}^N}\int_{\Omega}|x-\alpha|^2dx. 
\end{eqnarray*}
Then the lower bound (\ref{LiYau1}) can be improved as
\begin{eqnarray}
\label{Toni1}
&&\sum_{i=1}^m\lambda_{i}(\Omega)\geq\frac{NC_N}{N+2}\mu_N(\Omega)^{-\frac{2}{N}}m^{\frac{N+2}{N}}+M_N\frac{\mu_N(\Omega)}{I(\Omega)}\;m,\\
&&M_N=\frac{c}{N+2},\;\;\mbox{with $c<(2\pi)^2\omega_N^{-\frac{4}{N}}$, but $c$ independent of $N$}.\nonumber
\end{eqnarray}
 
In this paper we note that the improved estimate (\ref{Toni1}) implies a new and sharp estimate on the dimension of the global attractor, for the reaction diffusion system 
\begin{eqnarray}
\label{RC1}
\phi_t-\nu\Delta\phi-V(x)\phi&+&f(\phi)=0,\,\nu>0,\;\;\mbox{in}\;\;\Omega,\;
t>0,\\ 
\label{RC2}
\phi(x,0)&=&\phi_0(x),\;\;\mbox{for}\;\;x\in\Omega,\\
\label{RC2'}
\phi(x,t)&=&0\;\;\mbox{in}\;\;\partial\Omega,\;t>0,
\end{eqnarray}
where generally $\Omega\subset\mathbb{R}^N$ is a bounded domain of $N\geq 1$. The potential $V(x)$ is either {\em the inverse square potential} $V(x)=\delta\,|x|^{-2}$ or {\em the borderline potential} $V(x)=\delta\,\mathrm{dist}(x,\partial\Omega)^{-2}$, $\delta\geq 0$. For the reaction term $f:\mathbb{R}\rightarrow\mathbb{R}$ we are making the standard assumption that is a polynomial of odd degree with a positive leading coefficient, 
\begin{eqnarray}
\label{cnon}
f(s)=\sum_{k=1}^{2\gamma-1}b_ks^k,\;\;b_{2\gamma-1}>0.
\end{eqnarray}
With the nonlinearity (\ref{cnon}) equation (\ref{RC1}) can be considered as  a {\em singular Allen-Cahn-type equation} (including the classical Allen-Cahn equation when $\delta=0$). 
Setting $\mu:=\delta/\nu$, we consider as subcritical the cases $\delta=0$ for $N\geq1$ and $0<\mu<\mu^*$ for $N\geq 2$ or $N\geq 3$, depending on the type of the potential involved in (\ref{RC1}). The constant $\mu^*$ denotes the optimal constant of {\em Hardy and Hardy-type inequalities}. The critical value $\mu^*=(N-2)^2/4$ is the optimal constant of the Hardy inequality 
\begin{eqnarray}
\label{Har}
\mu\int_{\Omega}\frac{u^2}{|x|^2}dx\leq\int_{\Omega}|\nabla u|^2dx,\;\;\mbox{for all}\;\;u\in C^{\infty}_0(\Omega),\;\;\Omega\subseteq\mathbb{R}^N,\;N\geq 3,
\end{eqnarray}
which is not attained in $H^1_0(\Omega)$. Similarly, $\mu_*=1/4$ is the optimal constant of the Hardy-type inequality
\begin{eqnarray}
\label{Har2}
\mu\int_{\Omega}\frac{u^2}{d^2(x)}dx&\leq&\int_{\Omega}|\nabla u|^2dx,\;\;\mbox{for all}\;\;u\in C^{\infty}_0(\Omega),\;\;d(x):=\mathrm{dist}(x,\partial\Omega),\;\;\Omega\subset\mathbb{R}^N,\;N\geq 2,
\end{eqnarray}
which is also not attained in $H^1_0(\Omega)$.

This new estimate, although it simply comes out by incorporating (\ref{Toni1}) in the well known procedure estimating the distortion of infinitesimal $m$-volumes produced by the associated semiflow (cf. \cite{CFT1985,RTem88}), is sharp not only due to the optimality of the lower bound (\ref{Toni1}). Regarding its explicit dependence on the parameters of the problem and the geometric properties of the open set $\Omega$ {\em it reveals the existence of a threshold value on the ratio $\mathcal{R}(\Omega):=\mu_{N}(\Omega)/I(\Omega)$ on which the dimension of the attractor in the subcritical case may considerably change: \vspace{0.2cm}\newline 
There exists $\mathcal{R}_{\mathrm{thresh}}(N,\mu,f)>0$ such that if $\mathcal{R}(\Omega)<\mathcal{R}_{\mathrm{thresh}}$ then $\mathrm{dim}_{H}\mathcal{A}\leq d_0(N,\mu,f,\mathcal{R}(\Omega))$.  On the other hand, if $\mathcal{R}_{\mathrm{thresh}}\leq\mathcal{R}(\Omega)$ then  $\mathrm{dim}_{H}\mathcal{A}\leq 1$.}\vspace{0.2cm}\newline
The result is given in detail in Section \ref{mainSec}. In the case $\delta=0$, $N\geq 1$, the  result can be viewed as a new condition on the diffusivity $\nu$ for  a change of the attractor dimension (Theorem \ref{ACStand} \& Remark \ref{ACSa}). Furthermore it is verified that the dimension of the attractor should be actually smaller than the existing estimates indicate, see for example \cite{bab90,Chep99,RTem88}.

In Section \ref{ComSec}, we comment to the case of the critical potentials $\mu=\mu^*$. This case can be treated in generalized Sobolev spaces as the results of \cite{Ach03,brecab98,breduptes05,BrezisMarcus97,bv97,Ach07, vz00}, suggest.  Using the Weyl's type estimates on the eigenvalues of the critical Schr\"odinger operator $-\Delta-V$ derived in \cite{LMPNK}, we conclude by presenting the necessary conditions under which the dimension of the attractor in the critical case can be (explicitly) estimated.
\section{The subcritical case $0\leq \mu<\mu^*$: Sharp estimates on the Hausdorff dimension of the global attractor}
\label{mainSec}
\setcounter{equation}{0}
\paragraph{A. Inverse square potential $V(x)=\delta |x|^{-2}$.}
We consider first the case where $\Omega\subset\mathbb{R}^N$, $N\geq 3$, containing the origin and $\delta>0$ (the case $N=2$ reduces actually to the case $\delta=0$, see \cite[pg.108]{vz00}). Let us note that equation (\ref{RC1}) can be rewritten as
\begin{eqnarray*}
\partial_t\phi+\nu\mathcal{K}_s\phi+f(\phi)=0,\;\;\nu>0,
\end{eqnarray*}
where $\mathcal{K}_s$ denotes the Schr\"odinger operator
\begin{eqnarray*}
\mathcal{K}_s=-\Delta-\frac{\mu}{|x|^2},\;\;\mu=\frac{\delta}{\nu}.
\end{eqnarray*}
In the subcritical case $0<\mu<\mu^*=(N-2)^2/4$, the Hilbert space $H_{\mu}(\Omega)$ can be considered (cf.  \cite{vz00}), which is defined as the completion of $C^{\infty}_0(\Omega)$
with respect to the norm
\begin{eqnarray}
\label{PR2sub}
||u||_{H_{\mu}(\Omega)}:=\left[\int_{\Omega}\left(|\nabla u|^2-\frac{\mu}{|x|^2}u^2\right)dx\right]^{\frac{1}{2}},\;\;0<\mu<\mu^*.
\end{eqnarray}
Hardy's inequality (\ref{Har}) implies that the norm $||u||_{1,2}=\left(\int_{\Omega}|\nabla u|^2dx\right)^{1/2}$ of the Sobolev space $H^{1}_0(\Omega)$ and the norm (\ref{PR2sub}) are equivalent. One has that
\begin{eqnarray}
\label{PR3sub}
C_{\mu}\||u||^2_{H^{1}_{0}(\Omega)} \leq ||u||^2_{H_{\mu}(\Omega)} \leq
||u||^2_{H^{1}_{0}(\Omega)},
\end{eqnarray}
with $C_{\mu}: = 1 - \frac{\mu}{\mu^*}>0$, if $0<\mu<\mu^*$. The operator $\mathcal{K}_{s}:=-\Delta-\mu/|x|^2$ with domain 
\begin{eqnarray*}
D(\mathcal{K}_s):=\left\{u\in H_{\mu}(\Omega)\;:\;-\Delta u-\frac{\mu }{|x|^2}u\in L^2(\Omega),\;0<\mu<\mu^*\right\},
\end{eqnarray*}
is a nonnegative self-adjoint operator on $L^2(\Omega)$. Due to the standard embedding properties of $H^{1}_{0}(\Omega)\equiv H_{\mu}(\Omega)$ we have the following 
\begin{proposition}
\label{Subexist}
Let $\Omega\subset\mathbb{R}^N$, $N\geq 3$ be a bounded domain containing the origin, and assume that $0<\mu=\delta/\nu<\mu^*$. The problem (\ref{RC1})-(\ref{RC2'}) with nonlinearity (\ref{cnon}) and $V(x)=\delta |x|^2$, defines a semiflow $\mathcal{S}(t): L^2(\Omega)\rightarrow L^2(\Omega)$. The semiflow possesses a global attractor $\mathcal{A}_{\mu}$ which is bounded in $H^1_0(\Omega)$, compact and connected in $L^2(\Omega)$.  
\end{proposition}
We recall that since $f$ is of the form (\ref{cnon}), there exists $\kappa>0$ such that
\begin{eqnarray}
\label{derivcon}
f'(s)\geq -\kappa,\;\;\mbox{for all}\;s\in\mathbb{R}.
\end{eqnarray}
The main result result of the paper is
\setcounter{theorem}{1}
\begin{theorem}
\label{HDABD}
Let $\Omega\subset\mathbb{R}^N$, $N\geq 3$, a bounded domain containing the origin. Let $0<\mu=\delta/\nu<\mu^*$ and consider the global attractor $\mathcal{A}_{\mu}$ of the semiflow $\mathcal{S}(t):L^2(\Omega)\rightarrow L^2(\Omega)$. We define
\begin{eqnarray*}
\label{DefThresh}
\mathcal{R}_{\mathrm{thresh}}:=\frac{\kappa\mu^*}{M_N(\mu^*\nu-\delta)}.
\end{eqnarray*}
(i) If $\mathcal{R}_{\mathrm{thresh}}\leq \mathcal{R}(\Omega)=\mu_N(\Omega)/I(\Omega)$, 
then $\mathcal{A}_{\mu}$ has  finite Hausdorff dimension  $\mathrm{dim}_{H}\mathcal{A}_{\mu}\leq 1$.\newline
(ii) If $\mathcal{R}(\Omega)<\mathcal{R}_{\mathrm{thresh}}$,
then $\mathcal{A}_{\mu}$ has finite Hausdorff dimension $\mathrm{dim}_{H}\mathcal{A}_{\mu}\leq d_0$ and finite fractal dimension $\mathrm{dim}_{F}\mathcal{A}_{\mu}\leq d_0$ where  
\begin{eqnarray}
\label{HDA1''}
d_0= \left(\frac{N+2}{NC_N}\right)^{\frac{N}{2}}\left(\frac{\kappa\mu^*}{(\mu^*\nu-\delta)}-M_N\mathcal{R}(\Omega)\right)^{\frac{N}{2}}\mu_{N}(\Omega),\;\;N\geq 3.
\end{eqnarray}
\end{theorem}
{\bf Proof:} The usual arguments can be applied establishing that for every $t>0$,  the function $\phi_0\rightarrow\mathcal{S}(t)\phi_0$ is Fr\'{e}chet differentiable. The differential is $D(t,\phi_0):\xi\in L^2(\Omega)\rightarrow \Phi(t)\in L^2(\Omega)$, $\Phi(t)$  where $\Phi (t)$ is the solution of the first variation equation
\begin{eqnarray}
\label{fv1}
\Phi_t+\nu\mathcal{H}_s\Phi-f'(\mathcal{S}(t)\phi_0)\Phi=0\;\;\mbox{for}\;\;t>0,
\end{eqnarray}
supplemented with the initial and boundary conditions
\begin{eqnarray}
\label{fv2}
\Phi(t)=0\;\;\mbox{on}\;\;\partial\Omega\;\;\mbox{and}\;\;\Phi(0)=\xi\in L^2(\Omega).
\end{eqnarray} 
We shall examine the development of infinitesimal parallelepipeds spanned by $\Phi_1(t),\cdots\Phi_m(t)\in L^2(\Omega)$, where $\Phi_{i}(t),\;\;i=1,\ldots,m$ is an infinitesimal vector evolving from $\Phi_i(0)=\xi_i\in L^2(\Omega)$. The vectors $\Phi_{i}$, $i=1,\ldots,m$ are $m$-solutions of (\ref{fv1})-(\ref{fv2}) starting from the initial conditions $\Phi_i(0)=\xi_i\in L^2(\Omega)$ and it is well known that the $m$-dimensional volume $|\Phi_1(t)\wedge\Phi_2(t)\ldots\wedge\Phi_m(t)|$ of the infinitesimal parallelepiped spanned by $\Phi_i(t)$ is given by
\begin{eqnarray}
\label{fv3}
|\Phi_1(t)\wedge\Phi_2(t)\ldots\wedge\Phi_m(t)|&=&|\xi_1\wedge\ldots\wedge\xi_m|\exp\int_{0}^{t}\mathrm{Tr}[L(\mathcal{S}(s)\phi_0)\circ\mathcal{Q}_m(s)]ds,
\end{eqnarray}
where
\begin{eqnarray*}
L(\mathcal{S}(t)\phi_0)=-\nu\mathcal{H}_s\Phi-f'(\mathcal{S}(t)\phi_0)\Phi,\;\;\mbox{for}\;\;\Phi\in L^2(\Omega).
\end{eqnarray*}
We denote by $\mathcal{Q}_m(t)$, the orthogonal projection in $L^2(\Omega)$ onto $\mathrm{span}\left\{\Phi_1(t),\ldots,\Phi_m(t)\right\}$. We fix $t$ for the time being and we consider an orthonormal basis $e_1,e_2,\ldots$ of $L^2(\Omega)$ with $$\mathrm{span}\left\{e_1,\ldots,e_m\right\}=\mathrm{span}\left\{\Phi_1(t),\ldots,\Phi_m(t)\right\}=\mathrm{span}\mathcal{Q}_m(t)L^2(\Omega).$$
Since $\Phi_i(t)\in H_0^1(\Omega)$, for all $i\in\mathbb{N}$ and almost all $t>0$, we have $e_i\in H_0^1(\Omega)$ for all $i\in\mathbb{N}$. Then $\mathcal{Q}_m(t)e_{i}=e_{i}$ if $i\leq m$ and $\mathcal{Q}_m(t)e_{i}=0$ otherwise. Thus, for 
$i\leq m$,
\begin{eqnarray}
\label{fv4}
\left(L(\mathcal{S}(t)\phi_0)\circ\mathcal{Q}_m(t)e_i,e_i\right)_{L^2(\Omega)}&=&\left(L(\mathcal{S}(t)\phi_0)e_i,e_i\right)_{L^2(\Omega)}\nonumber\\
&=&-\nu\left\{\int_{\Omega}|\nabla e_i|^2dx-\int_{\Omega}\frac{\mu e^2_i}{|x|^2} dx\right\}-\int_{\Omega}f'(\mathcal{S}(t)\phi_0)e^2_i dx.
\end{eqnarray} 
Using (\ref{derivcon}), inequality (\ref{fv4}) becomes
\begin{eqnarray}
\label{fv5}
\left(L(\mathcal{S}(t)\phi_0)\circ\mathcal{Q}_m(t)e_i,e_i\right)_{L^2(\Omega)}
\leq-\nu\left\{\int_{\Omega}|\nabla e_i|^2dx-\int_{\Omega}\frac{\mu e^2_i}{|x|^2}dx\right\}+\kappa,
\end{eqnarray}
noticing that $\int_{\Omega}e_i^2dx=1$ by the orthonormality of the $e_i$. Therefore 
\begin{eqnarray}
\label{fv6'}
\mathrm{Tr}[L(\mathcal{S}(s)\phi_0)\circ\mathcal{Q}_m(s)]&=&\sum_{i=1}^m\left(L(\mathcal{S}(t)\phi_0)\circ\mathcal{Q}_m(t)e_i,e_i\right)_{L^2(\Omega)}\nonumber\\
&\leq&
-\nu\sum_{i=1}^{m}\left\{\int_{\Omega}|\nabla e_i|^2dx-\int_{\Omega}\frac{\mu e^2_i}{|x|^2} dx\right\}+\kappa m\nonumber\\
&=&-\nu\sum_{i=1}^m\left(\mathcal{K}_se_{i},e_{i}\right)_{L^2(\Omega)}+\kappa m\nonumber\\
&\leq&-\nu_1\sum_{i=1}^m\int_{\Omega}|\nabla e_i|^2dx+\kappa m,\;\;\nu_1:=\nu C_{\mu}\nonumber\\
&=&-\nu_1\sum_{i=1}^m(-\Delta e_i,e_i)_{L^2(\Omega)}+\kappa m.
\end{eqnarray}
Furthermore, by using the inequality
\begin{eqnarray*}
\sum_{i=1}^m\left(-\Delta e_{i},e_{i}\right)_{L^2(\Omega)}\geq\sum_{i=1}^m\lambda_i(\Omega),
\end{eqnarray*}
we may insert (\ref{Toni1}) in (\ref{fv6'}), to get
\begin{eqnarray}
\label{fv7'}
\mathrm{Tr}[L(\mathcal{S}(t)\phi_0)\circ\mathcal{Q}_m(t)]&\leq&
-\nu_1\frac{NC_N}{N+2}\mu_N(\Omega)^{-\frac{2}{N}}m^{\frac{N+2}{N}}-\nu_1 M_N\mathcal{R}(\Omega)\;m+\kappa\;m,\\
\label{fv7''}
&&=-\nu_1\frac{NC_N}{N+2}\mu_N(\Omega)^{-\frac{2}{N}}m^{\frac{N+2}{N}}+\kappa_1\;m,
\end{eqnarray}
where 
\begin{eqnarray}
\label{fv7a'}
\kappa_1=\kappa-\nu_1 M_N\mathcal{R}(\Omega).
\end{eqnarray}
Assume that $\kappa_1\leq 0$ i.e. $\mathcal{R}_{\mathrm{thresh}}\leq \mathcal{R}(\Omega)$.  Then (\ref{fv3}) implies the exponential decay in time of the $m$-dimensional volume $|\Phi_1(t)\wedge\Phi_2(t)\ldots\wedge\Phi_m(t)|$, {\em for any} $m\geq 1$. Then Constantin-Foias-Temam theory \cite[Proposition 2.1, pg. 364 \& Theorem 3.3, pg. 374]{RTem88} implies that the Hausdorff dimension of the attractor is less or equal than $m$, for any $m\geq 1$. This proves claim $(i)$ of the theorem.

Assume next that $\kappa_1>0$, i.e. condition $\mathcal{R}_{\mathrm{thresh}}> \mathcal{R}(\Omega)$ is satisfied. We  define by (\ref{fv7'}) the function
\begin{eqnarray*}
g(x)=-\beta\nu_1\mu_N(\Omega)^{-\frac{2}{N}}x^{\frac{N+2}{N}}+\kappa_1\;x,\;\;\beta=\frac{NC_N}{N+2}.
\end{eqnarray*}
This function has the root
\begin{eqnarray*}
d_0=\beta^{-\frac{N}{2}}\left(\frac{\kappa_1}{\nu_1}\right)^{\frac{N}{2}}\mu_N(\Omega).
\end{eqnarray*}
Then again the results of \cite{RTem88}, imply that $\mathrm{dim}_H\mathcal{A}_{\mu}\leq d_0$. In addition,  since the function $g$ is a concave function of the continuous variable $x$, \cite[Corollary 2.2, pg. 815]{Chep99} implies that the fractal dimension of $\mathcal{A}_{\mu}$ is $\mathrm{dim}_F\mathcal{A}_{\mu}\leq d_0$ 
This proves claim $(ii)$ of the theorem. \ $\blacksquare$

In the subcritical case $\delta=0$ (classical Allen-Cahn equation), we may consider the initial-boundary value problem (\ref{RC1})-(\ref{RC2'}) in a bounded domain $\Omega\subset\mathbb{R}^N$, $N\geq 1$. In this case we have
\begin{theorem}
\label{ACStand}
Let $\Omega\subset\mathbb{R}^N$, $N\geq 1$ be a bounded domain. Consider the global attractor $\mathcal{A}$ of the semiflow $\mathcal{S}(t):L^2(\Omega)\rightarrow L^2(\Omega)$ associated to (\ref{RC1})-(\ref{RC2'}) when $\delta=0$ and with the nonlinearity (\ref{cnon}). We define 
\begin{eqnarray*}
\label{DefThreshA}
\mathcal{R}_{\mathrm{thresh}}:=M_N\,\frac{\kappa}{\nu}.
\end{eqnarray*}
(i) Assume that $\mathcal{R}_{\mathrm{thresh}}\leq \mathcal{R}(\Omega)$. 
Then $\mathcal{A}$ has  finite Hausdorff dimension and $\mathrm{dim}_{H}\mathcal{A}\leq 1$.\newline
ii) Assume that $\mathcal{R}(\Omega)<\mathcal{R}_{\mathrm{thresh}}$.
Then $\mathcal{A}$ has finite Hausdorff dimension $\mathrm{dim}_{H}\mathcal{A}_{\mu}\leq d_0$ and finite fractal dimension $\mathrm{dim}_{F}\mathcal{A}_{\mu}\leq d_0$ where  
\begin{eqnarray}
\label{HDA1''A}
d_0= \left(\frac{N+2}{NC_N}\right)^{\frac{N}{2}}\left(\frac{\kappa}{\nu}-M_N\mathcal{R}(\Omega)\right)^{\frac{N}{2}}\mu_{N}(\Omega),\;\;N\geq 1.
\end{eqnarray}
\end{theorem}
\setcounter{remark}{3}
\begin{remark}
\label{ACSa}
We remark that the assumption $\mathcal{R}_{\mathrm{thresh}}\leq \mathcal{R}(\Omega)$ in Theorem \ref{ACStand} which can be viewed as
\begin{eqnarray}
\label{difsty}
\frac{M_N}{\mathcal{R}(\Omega)}\kappa\leq \nu,
\end{eqnarray}
gives a new condition on the diffusivity $\nu$ in order to have a global attractor of small Hausdorff dimension. We recall that if $\nu$ is sufficiently large (cf. \cite[Remark 1.2, pg. 88]{RTem88}), $||\phi(t)||_{L^2(\Omega)}^2$ decreases exponentially to $0$ as $t\rightarrow\infty$, without necessarily implying that $\mathcal{A}=\{0\}$. The attractor $\mathcal{A}$ may contain one or many heteroclinic curves. 
On the other hand, when  $\mathcal{R}(\Omega)<\mathcal{R}_{\mathrm{thresh}}$ meaning that
\begin{eqnarray}
\label{difstyb}
\frac{M_N}{\mathcal{R}(\Omega)}\kappa>\nu,
\end{eqnarray}
we get an improved upper bound (\ref{HDA1''}) compared with  existing upper bounds on the dimension of the global attractor for the Allen-Cahn equation \cite{Chep99, RTem88}. The estimate (\ref{HDA1''A}) shows that the Hausdorff dimension is indeed smaller due to the appearance of the term $-M_N\mathcal{R}(\Omega)$. For the case $\delta>0$ the conditions of Theorem  \ref{HDABD}, can be similarly implemented.
\end{remark}
\paragraph{B. Borderline potential $V(x)=\delta d^{-2}(x)$.}
In the case of the borderline potential, we start with the case of $\Omega\subset\mathbb{R}^N$, $N\geq 2$. Now equation (\ref{RC1}) can be rewritten as
\begin{eqnarray*}
\partial_t\phi+\nu\mathcal{H}_s\phi+f(\phi)=0,\;\;\nu>0,
\end{eqnarray*}
where $\mathcal{K}_s$ denotes the Schr\"odinger operator
\begin{eqnarray*}
\mathcal{H}_s=-\Delta-\frac{\mu}{d^2(x)},\;\;\mu=\frac{\delta}{\nu}.
\end{eqnarray*}
In the subcritical case $0<\mu<\mu^*=1/4$, and motivated by \cite{BrezisMarcus97}, we may consider the Hilbert space $W_{\mu}(\Omega)$  defined as the completion of $C^{\infty}_0(\Omega)$
with respect to the norm
\begin{eqnarray}
\label{PR2sub''}
||u||_{W_{\mu}(\Omega)}:=\left[\int_{\Omega}\left(|\nabla u|^2-\mu\frac{u^2}{d^2(x)}\right)dx\right]^{\frac{1}{2}},\;\;0<\mu<\mu^*.
\end{eqnarray}
Similarly to the inverse potential case, the Hardy-type inequality (\ref{Har2}) implies that the usual Sobolev norm and the norm (\ref{PR2sub''}) are equivalent. One has that
\begin{eqnarray*}
C_{\mu}\||u||^2_{H^{1}_{0}(\Omega)} \leq ||u||^2_{W_{\mu}(\Omega)} \leq
||u||^2_{H^{1}_{0}(\Omega)},
\end{eqnarray*}
this time with $C_{\mu}: = 1 - 4\mu> 0$, if $0<\mu<\mu^*$. The operator $\mathcal{H}_{s}:=-\Delta-\frac{\mu}{d^2(x)}$ with domain 
\begin{eqnarray*}
D(\mathcal{H}_s):=\left\{u\in W_{\mu}(\Omega)\;:\;-\Delta u-\mu\frac{u}{d^2(x)}\in L^2(\Omega),\;0<\mu<\mu^*\right\},
\end{eqnarray*}
is a nonnegative self-adjoint operator on $L^2(\Omega)$. Working exactly as for the proof of Theorem \ref{HDABD} we have
\setcounter{theorem}{4}
\begin{theorem}
\label{HBord}
Let $\Omega\subset\mathbb{R}^N$, $N\geq 2$, a bounded domain. Let $0<\mu=\delta/\nu<\mu^*$ and consider the global attractor $\mathcal{A}_{\mu}$ of the semiflow $\mathcal{S}(t):L^2(\Omega)\rightarrow L^2(\Omega)$ associated to  (\ref{RC1})-(\ref{RC2'}) with the nonlinearity (\ref{cnon}) and $V(x)=\delta d^{-2}(x)$. We define 
\begin{eqnarray*}
\label{DefThreshBP}
\mathcal{R}_{\mathrm{thresh}}:=\frac{\kappa}{M_N(\nu-4\delta)}.
\end{eqnarray*}
(i) Assume that $\mathcal{R}_{\mathrm{thresh}}\leq \mathcal{R}(\Omega)$. 
Then $\mathcal{A}_{\mu}$ has  finite Hausdorff dimension  $\mathrm{dim}_{H}\mathcal{A}_{\mu}\leq 1$.\newline
(ii) Assume that $\mathcal{R}(\Omega)<\mathcal{R}_{\mathrm{thresh}}$.
Then $\mathcal{A}_{\mu}$ has finite Hausdorff dimension $\mathrm{dim}_{H}\mathcal{A}_{\mu}\leq d_0$ and finite fractal dimension $\mathrm{dim}_{F}\mathcal{A}_{\mu}\leq d_0$ where  
\begin{eqnarray*}
\label{HDA1''BP}
d_0= \left(\frac{N+2}{NC_N}\right)^{\frac{N}{2}}\left(\frac{\kappa}{(\nu-4\delta)}-M_N\mathcal{R}(\Omega)\right)^{\frac{N}{2}}\mu_{N}(\Omega),\;\;N\geq 2.
\end{eqnarray*}
\end{theorem}
\section{Remarks on the critical cases $\mu=\mu^*$}
\label{ComSec}
\setcounter{equation}{0}
In the critical case $\mu=\mu^*$ the initial-boundary value problem (\ref{RC1})-(\ref{RC2'}) still defines a semiflow $\mathcal{S}(t):L^2(\Omega)\rightarrow L^2(\Omega)$. For the definition of the semiflow and the existence and finite dimensionality of the global attractor generalized Sobolev spaces come into play, as well as,  Weyl's type estimates on the eigenvalues of the critical Schr\"odinger operators. These results can be used under some further geometric restrictions on $\Omega$. 
\paragraph{A. Critical inverse square potential $V(x)=\frac{(N-2)^2}{4|x|^2}$.}
We assume that $\Omega\subset\mathbb{R}^N$, $N\geq 3$, is a bounded domain containing the origin. In the case of the inverse square potential, the critical initial-boundary value problem is for the equation
\begin{eqnarray}
\label{cris}
\partial_t\phi+\nu\mathcal{K}\phi+f(\phi)=0,\;\;\nu>0,
\end{eqnarray}
supplemented with the initial and boundary conditions (\ref{RC2})-(\ref{RC2'}). Here $\mathcal{K}$ denotes the critical Schr\"odinger operator $-\Delta-\frac{(N-2)^2}{4|x|^2}$. The operator $\mathcal{K}:D(\mathcal{K})\rightarrow L^2(\Omega)$, with its domain defined as 
\begin{eqnarray*}
\label{cris1}
D(\mathcal{K}):=\left\{u\in H(\Omega)\;:\;-\Delta u-\frac{(N-2)^2}{4|x|^2}u\in L^2(\Omega)\right\},
\end{eqnarray*}
is a nonnegative self-adjoint operator on $L^2(\Omega)$. Here $H(\Omega)$ is the Hilbert space defined in \cite[Section 4.1]{vz00}, as the completion of $C^{\infty}_0(\Omega)$ in the norm 
\begin{eqnarray}
\label{PR2}
||u||_{H(\Omega)}:=\left[\int_{\Omega}\left(|\nabla u|^2-\frac{(N-2)^2}{4|x|^2}u^2\right)dx\right]^{\frac{1}{2}}. 
\end{eqnarray}
The following crucial improvement of (\ref{Har}) has been proved in \cite[Theorem 2.2, pg. 108]{vz00}: There exists some positive constant $C(r,\Omega)$ such that
\begin{eqnarray}
\label{PR4}
\int_{\Omega}|\nabla u|^2dx-\int_{\Omega}\frac{(N-2)^2}{4|x|^2}u^2dx\geq C(r,\Omega) ||u||^2_{W^{1,r}(\Omega)},\;\;\mbox{for all}\;\;u\in C^{\infty}_0(\Omega),\;\;1\leq r<2,
\end{eqnarray}
implying that the embedding $H(\Omega)\hookrightarrow W^{1,r}(\Omega),\;\;1\leq r<2$. From this, we actually infer the compact embeddings $H(\Omega) \hookrightarrow\hookrightarrow
L^2(\Omega)$, $H(\Omega) \hookrightarrow\hookrightarrow
H_{0}^{s}(\Omega),\,0\leq s<1$, and hence
the existence of a complete orthonormal basis $\left\{\phi_j\right\}_{j\geq 1}$ of $L^2(\Omega)$ consisting of eigenfunctions of $\mathcal{K}$ with the eigenvalue sequence
\begin{eqnarray}
\label{PR6}
0<\lambda_1\leq\lambda_2\leq\cdots\leq\lambda_j\leq\cdots\rightarrow\infty,\;\;\mbox{as}\;\;j\rightarrow\infty.
\end{eqnarray}
Furthermore, it was shown in \cite[Theorem 2.1]{LMPNK} that if
\begin{eqnarray}
\label{crExpA}
\frac{2N}{N+2}< q <2,
\end{eqnarray} 
the eigenvalues (\ref{PR6}) satisfy the Weyl's estimate 
\begin{eqnarray}
\label{Weyl's1}
\lambda_j\geq C(q,\Omega)\mathrm{e}^{-1}\mu_N(\Omega)^{-\frac{Nq-2N+2q}{Nq}}j^{\frac{Nq-2N+2q}{Nq}},\;\;j\rightarrow\infty.
\end{eqnarray}
With the estimate (\ref{Weyl's1}) in hand we have
\begin{theorem}
Let $\Omega\subset\mathbb{R}^N$, $N\geq 3$ a bounded domain containing the origin. The initial-boundary value problem (\ref{cris1})-(\ref{RC2})-(\ref{RC2'}) with nonlinearity (\ref{cnon}), defines a semiflow $\tilde{\mathcal{S}}(t): L^2(\Omega)\rightarrow L^2(\Omega)$ possessing a global attractor $\tilde{\mathcal{A}}$.  There exists a constant $C_1(q,\Omega, N)>0$ such that $\mathrm{dim}_H\tilde{\mathcal{A}}\leq \tilde{d}_0$ with 
\begin{eqnarray}
\label{dimis}
\tilde{d}_0=C_1(q,\Omega, N)^{-\frac{Nq}{Nq-2N+2q}}\left(\frac{\kappa}{\nu}\right)^{\frac{Nq}{Nq-2N+2q}}\mu_N(\Omega).
\end{eqnarray}
\end{theorem}
{\bf Proof:} We just note that the process followed in Theorem \ref{HDABD}, leads to the corresponding estimate for (\ref{cris})
\begin{eqnarray}
\label{cris3}
\mathrm{Tr}[L(\mathcal{S}(t)\phi_0)\circ\mathcal{Q}_m(t)]\leq -\frac{\nu C_1(q,\Omega)}{\mu_N(\Omega)^{\frac{p^*-2}{p^*}}}\;m^{\frac{2(p^*-1)}{p^*}}+\kappa m,\;\;p^*=\frac{qN}{N-q}.
\end{eqnarray}
We remark that the requirement $\frac{2(p^*-1)}{p^*}>1$ justifies  condition (\ref{crExpA}). The right-hand side of (\ref{cris3}) has the root $\tilde{d}_0$ and the estimate on the dimension follows from \cite{RTem88}, \cite[Corollary 2.2, pg. 815]{Chep99}.\ \ $\blacksquare$

Note that the estimate (\ref{dimis}) lacks the exponent $N/2$ since the improved Hardy-Poincar\'{e} inequality (\ref{PR4}) is not valid in the critical value $r=2$ and consequently, Weyl's estimate (\ref{Weyl's1}) is not valid in the critical value $q=2$. This is not the case for the subcritical problem $0<\mu<\mu^*$.
\paragraph{B. Critical  borderline potential $V(x)=\frac{1}{4d^2(x)}$.}
We assume in this case that $\Omega\subset\mathbb{R}^N$, $N\geq2$,  is a bounded  \emph{ smooth and convex} domain. In the case of the critical borderline potential, the initial-boundary value problem is for the equation
\begin{eqnarray}
\label{crisbp}
\partial_t\phi+\nu\mathcal{H}\phi+f(\phi)=0,\;\;\nu>0,
\end{eqnarray}
supplemented with the initial and boundary conditions (\ref{RC2})-(\ref{RC2'}). Now $\mathcal{H}$ denotes the critical the critical Schr\"odinger operator $-\Delta-\frac{1}{4d^2(x)}$. The operator $\mathcal{H}:D(\mathcal{H})\rightarrow L^2(\Omega)$, with  domain 
\begin{eqnarray*}
\label{PR3'}
D(\mathcal{H}):=\left\{u\in W(\Omega)\;:\;-\Delta u-\frac{u}{4d^2(x)}\in L^2(\Omega)\right\},
\end{eqnarray*}
is a nonnegative self-adjoint operator on $L^2(\Omega)$. The Hilbert space $W(\Omega)$ is the  completion of $C^{\infty}_0(\Omega)$  in the norm
\begin{eqnarray*}
\label{PR2'}
||u||_{W(\Omega)}:=\left[\int_{\Omega}\left(|\nabla u|^2-\frac{u^2}{4d^2(x)}\right)dx\right]^{\frac{1}{2}}.
\end{eqnarray*}
Improving the results of \cite{BrezisMarcus97}, it was shown in \cite{Ach03}, that there exists some positive constant $C$ such that
\begin{eqnarray*}
\int_{\Omega}|\nabla u|^2dx-\frac{1}{4}\int_{\Omega}\frac{u^2}{d^2(x)}
dx\geq C\left(\int_{\Omega}|\nabla u|^q dx\right)^{\frac{2}{q}},\;\;\mbox{for all}\;\;u\in C^{\infty}_0(\Omega),\;\;1\leq q<2.
\end{eqnarray*}
Thus as in A. we have the compact embeddings $H(\Omega) \hookrightarrow\hookrightarrow
L^2(\Omega)$, $H(\Omega) \hookrightarrow\hookrightarrow
H_{0}^{s}(\Omega),\,0\leq s<1$, and the existence of a complete orthonormal basis $\{\tilde{\phi}_j\}_{j\geq 1}$ of $L^2(\Omega)$ made of eigenfunctions of $\mathcal{H}$ with the eigenvalue sequence
\begin{eqnarray}
\label{PR6'}
0<\tilde{\lambda}_1\leq\tilde{\lambda}_2\leq\cdots\leq\tilde{\lambda}_j\leq\cdots\rightarrow\infty,\;\;\mbox{as}\;\;j\rightarrow\infty.
\end{eqnarray}
The {\em convexity condition can not be relaxed} since it is known that if $\Omega$ is not convex we may have $\lambda_1>-\infty$, \cite{BrezisMarcus97}. 
By making use of the inequality in \cite[Theorem 1.1, pg. 492]{Ach06}, 
\begin{eqnarray}
\label{PR7'}
\int_{\Omega}|\nabla u|^pdx-\left(\frac{p-1}{p}\right)^p\int_{\Omega}\frac{|u|^p}{d^p}dx\geq \tilde{C}(\Omega)\left(\int_{\Omega}|u|^rdx\right)^{\frac{p}{r}},\;\;\mbox{for all}\;\;u\in C^{\infty}_0(\Omega),
\end{eqnarray}
which holds for $1<p<N$ and $p\leq r<\frac{Np}{N-p}$ and its sharp estimates on the optimal constant 
\begin{eqnarray}
\label{PR8'}
c_1(p,r,N)D^{n-p-\frac{Np}{r}}_{\mathrm{int}}\geq \tilde{C}(\Omega)\geq c_2(p,r,N)D^{n-p-\frac{Np}{r}}_{\mathrm{int}},\;\;D_{\mathrm{int}}:=2\sup_{x\in\Omega}d(x),
\end{eqnarray}
it was shown in \cite[Theorem 2.3]{LMPNK} that if (\ref{crExpA}) holds, the eigenvalues (\ref{PR6'}) satisfy the Weyl's type estimate
\begin{eqnarray}
\label{Weyl's2'}
\tilde{\lambda}_j\geq \hat{C}(\Omega)\mathrm{e}^{-1}\mu_N(\Omega)^{-\frac{Nq-2N+2q}{Nq}}j^{\frac{Nq-2N+2q}{Nq}},\;\;j\rightarrow\infty.
\end{eqnarray}
The constant $\hat{C}(\Omega)$ satisfies the upper and lower estimates
\begin{eqnarray}
\label{Weyl's2''}
c_1(q,N)D^{\frac{N(q-2)}{q}}_{\mathrm{int}}\geq \hat{C}(\Omega)\geq c_2(q,N)D^{\frac{N(q-2)}{q}}_{\mathrm{int}}.
\end{eqnarray}
However, when $N\geq 3$ due to a further improvement of (\ref{PR8'}) to the critical value $q=2$, given in \cite[Theorem 3.4, pg. 46]{ACH007} ,
\begin{eqnarray}
\label{CrA1}
\int_{\Omega}|\nabla u|^2dx-\frac{1}{4}\int_{\Omega}\frac{u^2}{d^2}dx\geq C(N,D_{\mathrm{int}})\left(\int_{\Omega}|u|^{\frac{2N}{N-2}}dx\right)^{\frac{N-2}{N}},\;\;\mbox{for all}\;\;u\in C^{\infty}_0(\Omega),
\end{eqnarray} 
we get the Weyl's type estimate \cite[Theorem 2.4]{LMPNK}
\begin{eqnarray}
\label{CritW}
\tilde{\lambda}_j\geq C(N, D_{\mathrm{int}})\mathrm{e}^{-1}\mu_N(\Omega)^{-\frac{2}{N}}j^{\frac{2}{N}},\;\;j\rightarrow\infty.
\end{eqnarray}
With the estimates (\ref{Weyl's2'}) and (\ref{CritW}) we have 
\begin{theorem}
(i) Let $\Omega\subset\mathbb{R}^N$, $N\geq 2$ a bounded, smooth and convex domain. The initial-boundary value problem (\ref{crisbp})-(\ref{RC2})-(\ref{RC2'}) with nonlinearity (\ref{cnon}), defines a semiflow $\hat{\mathcal{S}}(t): L^2(\Omega)\rightarrow L^2(\Omega)$ possessing a global attractor $\hat{\mathcal{A}}$.  There exists a constant $C_1(q, N,D_{\mathrm{int}})>0$ such that $\mathrm{dim}_H\tilde{\mathcal{A}}\leq \hat{d}_0$ with 
\begin{eqnarray}
\hat{d}_0=C_1(q,N,D_{\mathrm{int}} )^{-\frac{Nq}{Nq-2N+2q}}\left(\frac{\kappa}{\nu}\right)^{\frac{Nq}{Nq-2N+2q}}\mu_N(\Omega).
\end{eqnarray}
(ii) Let $\Omega\subset\mathbb{R}^N$, $N\geq 3$ a bounded, smooth and convex domain. Then there exists $C_1(N,D_{\mathrm{int}})>0$ such that
$\mathrm{dim}_H\tilde{\mathcal{A}}\leq \hat{d}_0$ and  $\mathrm{dim}_F\tilde{\mathcal{A}}\leq \hat{d}_0$ with 
\begin{eqnarray}
\hat{d}_0=C_1(N,D_{\mathrm{int}} )^{-\frac{N}{2}}\left(\frac{\kappa}{\nu}\right)^{\frac{N}{2}}\mu_N(\Omega).
\end{eqnarray}
\end{theorem}
We conclude by mentioning that in the critical cases $\mu=\mu^*$, a semiflow may not be defined in $H^1_0(\Omega)$ as the non-existence results of \cite{brecab98,bv97,vz00} indicate. However, the semiflows can be defined in the corresponding generalized Sobolev phase spaces $H(\Omega),W(\Omega)$ (possibly under appropriate smallness conditions on the growth of the nonlinearity) and they satisfy the appropriate energy equations \cite{Ball1a,Ball2}. We also refer to our recent work \cite{NKNZ08}. 
\bibliographystyle{amsplain}

\end{document}